\documentclass[12pt]{amsart}
\usepackage{amssymb,xspace,graphicx,enumerate}
\usepackage[author-year]{amsrefs}
\newtheorem{mainthm}{Theorem}
\newtheorem{thm}{Theorem}[section]
\newtheorem{prop}[thm]{Proposition}
\newtheorem{lem}[thm]{Lemma}

\theoremstyle{definition}

\theoremstyle{remark}
\newtheorem{rmk}{Remark}[section]

\newtheorem*{rmk*}{Remark}
\newtheorem*{rmks*}{Remarks}
\newenvironment{enumeratei}{\begin{enumerate}[\upshape (i)]}
                           {\end{enumerate}}

\newenvironment{enumeraten}{\begin{enumerate}[{\tt 1.}]}
                           {\end{enumerate}}
\newcommand{\deq}{\overset{{\rm def}}{=}}            % defined equal
\newcommand{\E}{E}                                   % expectation
\renewcommand{\P}{P}                                 % probability
\newcommand{\Pb}[1]{\P\!\left( #1 \right)}           %   with brackets

\newcommand{\N}{{\mathbb N}}                         % natural numbers

\newcommand{\R}{{\mathbb R}}                         % real numbers
\newcommand{\F}{{\mathcal F}}                        % sigma-field
\newcommand{\Given}{\; \vline \;}

\newcommand{\one}{{\bf 1}}                           % indicator function

\renewcommand{\and}{\; \mbox{ and } \;}

\renewcommand{\emptyset}{\varnothing}
\newcommand{\ep}{\epsilon}

\renewcommand{\phi}{\varphi}
\newcommand{\convfdd}{\stackrel{{\rm f.d.d.}}{\longrightarrow}}
\newcommand{\weakconv}{\stackrel{{\rm d}}{\longrightarrow}}

\begin{document}
   \title[Hypergraph Processes]{Continuous And Discontinuous Phase 
    Transitions In Hypergraph Processes}

   \author[R.W.R.\@ Darling]{R.W.R.\@ Darling}
   \address{National\@ Security\@ Agency\\
        P.O.\ Box 535\\ Annapolis Junction, MD 20701}
   \email{rwrd@afterlife.ncsc.mil}
   %\urladdr{}
%
   \author[D.A.\@ Levin]{David A.\@ Levin}
   \address{Department of Mathematics\\The University\@ of Utah\\
      155 S.\@ 1400 E.\\Salt Lake City, UT 84112--0090}
   \email{levin@math.utah.edu}
   \author[J.R.\@ Norris]{James R.\@ Norris}
   \address{Statistical Laboratory\\ Centre For Mathematical Sciences\\
     Wilberforce Road, Cambridge, CB3 0WB}
   \email{j.r.norris@statslab.cam.ac.uk}

\maketitle

\vspace{-0.25in}

\centerline{\emph{National Security Agency, University of Utah, 
University of Cambridge}}

\begin{abstract}
  Let $V$ denote a set of $N$ vertices. To construct a \emph{hypergraph
  process}, create a new hyperedge at each event time of a Poisson
  process; the cardinality $K$ of this hyperedge is random, with
  generating function $\rho(x) \deq \sum \rho_k x^k$,
  where $\Pb{K = k} = \rho_k$; given $K=k$, the $k$ vertices appearing 
  in the new hyperedge are selected uniformly at random from
  $V$. Assume $\rho_1 + \rho_2 > 0$. Hyperedges of cardinality $1$ 
  are called \emph{patches}, and serve as a way of selecting root
  vertices. Identifiable vertices are those which are reachable
  from these root vertices, in a strong sense which generalizes the 
  notion of graph component. Hyperedges are called identifiable if all 
  of their vertices are identifiable. We use ``fluid limit'' scaling: 
  hyperedges arrive at rate $N$, and we study structures of size 
  $O(1)$ and $O(N)$. After division by $N$, numbers of identifiable 
  vertices and hyperedges exhibit phase transitions, which 
  may be continuous or discontinuous depending on the shape of the
  structure function $-\log(1 - x)/\rho'(x), \ x\in(0,1)$. 
  Both the case $\rho_1 > 0$, and the case $\rho_1 = 0 < \rho_2$ 
  are considered; for the latter, a single extraneous patch is added
  to mark the root vertex. 
\end{abstract}

\newpage
%\tableofcontents

\section{Introduction}

The \emph{$k$-core} of a graph is the largest subgraph with
minimum degree at least $k$.
\ocite{PSW:core} study the following algorithm for
finding the $2$-core of a graph:
\begin{enumeraten}
\item If vertices of degree one exist, select one
  and remove the edge incident to it.  This may
    cause the degree of other vertices to drop.
\item If there are no degree one vertices remaining,
  stop.
\item Repeat.
\end{enumeraten}
The graph obtained at the conclusion of this
algorithm is the $2$-core. 

This algorithm is a special case of another,
run on hypergraphs, called \emph{hypergraph collapse} and first
studied in \ocite{DN:HG}. 
By a \emph{hypergraph} we shall mean a map $\Lambda:2^V\to\{0,1,2,\dots\}$, 
where $V$ is a finite set of \emph{vertices} and $2^V$ is the set of subsets 
of $V$.
It will sometimes be helpful to think in terms of an {\em edge-labelling} of
$\Lambda$, which is a choice of a set $I$ and a map 
$e:I\to 2^V$ such that $\Lambda(A)=|\{i\in I: e(i)=A\}|$ for all $A$.
Thus $e$ describes a set of labelled subsets of $V$, which we call
\emph{hyperedges} and then $\Lambda$ gives
the number of hyperedges at each subset of $V$. 
Hyperedges of unit cardinality are called \emph{patches}.  Hypergraph
collapse is the following algorithm: 
\begin{enumeraten}
\item If a patch exists, select one and remove it together with the unique 
  vertex $v$ it contains. This will cause
  any other hyperedge $e(i)$ containing $v$ to be replaced by
  $e(i) \setminus \{v\}$.
\item If there are no patches remaining, stop.
\item Repeat.
\end{enumeraten}
Although we have described the algorithm in terms of an edge-labelled
hypergraph, 
the possible moves for $\Lambda$ do not depend on the edge-labelling chosen.
The vertices which are removed by hypergraph collapse are called
\emph{identifiable}, and hyperedges which contain only identifiable
vertices are also called \emph{identifiable}. These definitions do not depend
on the order in which patches are chosen during hypergraph collapse;
see \ocite{DN:HG}.  

The core-finding algorithm of \ocite{PSW:core} is hypergraph collapse
applied to the \emph{dual} hypergraph. To obtain the dual, note that we can
think of $e$ as a subset of $V\times I$. The roles of $V$ and $I$ are now 
symmetric, so $e$ also corresponds to an edge-labelling of a hypergraph
$\Lambda'$ 
in which the status of vertices and hyperedges is reversed.
Vertices (resp.\ hyperedges) of
$\Lambda$ not in the core correspond to identifiable hyperedges (resp.\
identifiable vertices) of $\Lambda'$. More information about graph
cores can be found in \ocite{F:core}, and hypergraph cores are
considered by \ocite{C:core}.  

The identifiable vertices obtained by hypergraph collapse
also serve to generalize to hypergraphs the definition of
graph component. A graph is a hypergraph having edges only of
cardinality two, and consequently has no patches.  However, if the
single hyperedge $\{v\}$ is added to the graph [making it a
  hypergraph], then the identifiable vertices obtained by running
hypergraph collapse on the augmented graph are exactly the vertices in
the graph component containing $v$. The identifiable edges are all the
edges of this graph component.   

This motivates the following definition for patch-free hypergraphs:  A
vertex is \emph{in the domain} 
of $v$ if it is in the set of identifiable vertices when
the hypergraph is augmented by the addition of the
hyperedge $\{v\}$.  

The purpose of this paper is to study the time-evolution of the
set of identifiable vertices and the set of identifiable edges
in a \emph{Poisson hypergraph process}, which is a hypergraph-valued,
continuous-time stochastic process.   The vertex set is $V =
\{1,2,\ldots,N\}$, and the process depends on parameters
$\{\rho_j\}_{j=1}^N$. Attached to each subset $A$ of $V$ is a Poisson
clock run at rate $N \rho_{|A|}/ \binom{N}{|A|}$, and these
clocks are independent of one another. [Here $|A|$ denotes
the cardinality of $A$.]  When the clock associated to $A$ ``rings'',
a new hyperedge equal to $A$ is added to the hypergraph.  The overall rate
at which hyperedges of cardinality $k$ are added is then $N \rho_k$.
We will call this process the Poisson($\rho$) hypergraph process.
This is a generalization of the ordinary random graph process, in
which edges form between each pair of vertices independently at a fixed rate.

While for $N$ fixed, this process depends only on the finite sequence
$\{\rho_k\}_{k=1}^N$, we will be interested in the asymptotic behavior
as $N \rightarrow \infty$, so we will assume that always the infinite
sequence $\{\rho_k\}_{k=1}^\infty$ is given.  Moreover, this sequence
is required to be a probability distribution on $\{1,2,\ldots\}$ with
finite expectation and satisfying $\rho_1 + \rho_2 > 0$.  The
generating function $x \mapsto \sum_{k=1}^\infty \rho_k x^k$ will
be denoted by $\rho$.  

In \ocite{DN:HG}, the \emph{Poisson($\beta$) random hypergraph} is
defined, where $\{\beta_k\}$ is a sequence of positive
real numbers.  This is a random hypergraph with vertex set
V = $\{1,2,\ldots,N\}$, so that for each $A \subset V$,
the number of occurrences of the hyperedge $A$ is a
Poisson random variable with expectation $N\beta_{|A|}/\binom{N}{|A|}$,
and these random variables are independent for different subsets of $V$.
If $\{\Lambda_t\}_{t \geq 0}$ is a Poisson($\rho$) hypergraph
process, then for fixed $t \geq 0$, $\Lambda_t$ is a Poisson($t\rho$)
random hypergraph.

We separate out two distinct cases in the study of Poisson hypergraph
processes, depending on whether $\rho_1 > 0$ or $\rho_1 = 0$. When
$\rho_1 = 0$, the hypergraph never acquires patches, and provided the
initial hypergraph is patch-free, the set of identifiable vertices
is forever void. As the previous discussion of ordinary graphs
suggests,  it is natural to consider in such cases the set of vertices
in the domain of a distinguished vertex. 

We discuss now the case $\rho_1 > 0$.
Our first result describes the evolution of the rescaled number of
identifiable vertices and hyperedges in the Poisson($\rho$)
hypergraph process $\{\Lambda_t\}_{t \geq 0}$.  Let 
\begin{equation} \label{eq:T_tilde_N_and_Z_tilde_N_defn}
  \begin{split}
    \tilde{T}^N_t & = \frac{|\text{identifiable vertices in } \Lambda_t|}{N}\\
    \tilde{Z}^N_t & = \frac{|\text{identifiable hyperedges in } \Lambda_t|}{N} \,.
  \end{split}
\end{equation}
The \emph{structure function} $t$, defined as
\begin{equation} \label{eq:structure_function}
  t(x) \deq \frac{ - \log(1-x) } { \rho'(x) } \,, \quad
  x \in (0,1) \,,
\end{equation}
plays a central role for hypergraph processes.
[Recall that $\rho(x) = \sum_{k=1}^\infty \rho_k x^k$.]
Typically $t$ is not invertible, but there
is a right-continuous monotonic function called the 
\emph{lower envelope}:
\begin{equation} \label{eq:gsdefinition}
  g(s) \deq \inf\{ x \in (0,1) \;:\; t(x) > s\}\,, 
  \quad s \geq 0 \,.
\end{equation}
Also important for hypergraph processes is 
the \emph{upper envelope}: 
\begin{equation} \label{eq:g_star_definition}
g^\star(s) \deq 
    \sup\{x \in (0, 1)  \;:\; t(x) < s\} \vee 0 \,, 
    \quad s \geq 0 \,.
\end{equation}

\begin{figure}
  \begin{center}
    \includegraphics{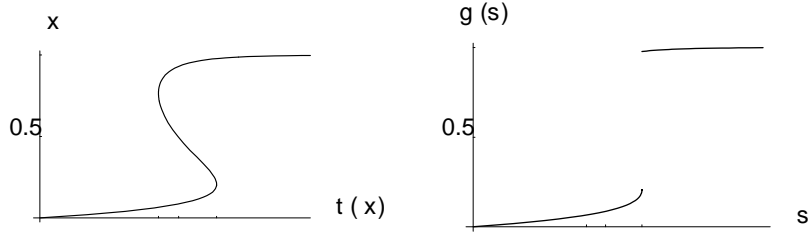}
    \caption{\emph{Left}: 
      Bicritical structure function, with $t(x)$ on the horizontal axis,
      corresponding to a quartic polynomial $\rho(x)$ with 
      $0 < \rho_1 < \rho_2 < \rho_3 < \rho_4$. \emph{Right}: 
      Lower envelope, showing the single discontinuity. 
      \label{fig:bicritical}}
  \end{center}
\end{figure}

We classify structure functions into three types:
graph-like, bicritical, and exceptional.  This
taxonomy is given in Table \ref{tab:class}. Figure
\ref{fig:bicritical} shows a bicritical structure function and the
corresponding lower envelope.  

\begin{table}[h]
  \begin{center}
    \begin{tabular}{c|c|c}
      \textsc{Type} & \textsc{Description} & \textsc{Example of $\rho(x)$} \\
      \hline
      graph-like & \parbox{2.00in}{\vspace{0.1in}
	$t$ is strictly increasing, and \\
	$g$ and $g^\star$ are continuous.\vspace{0.1in}}
      & cubic with $3 \rho_3 \leq \rho_2$ \\
      \hline
      bicritical & \parbox{2.00in}{\vspace{0.1in}
	$g$ and $g^\star$ each have \\ 
	exactly one discontinuity.
	\vspace{0.1in} }
      & cubic with $3 \rho_3 > \rho_2$ \\
      \hline
      exceptional & \parbox{2.00in}{\vspace{0.1in}
	$g$ or $g^\star$ has two or more \\
	discontinuities.
	\vspace{0.1in} }
      & $\frac{x + 5x^3 + 994x^{200}}{1000}$ 
    \end{tabular}
    
    \vspace{0.1in}

    \caption{Classification of structure functions \label{tab:class}}
  \end{center}
\end{table}

Let $\Xi \subset \R_+$ denote the discontinuity set of
$g$:
\begin{equation} \label{eq:discontinuity_set}
  \Xi \deq \{ s > 0 \,:\, g(s-) \neq g(s) \} \,,
\end{equation}
where  $g(s-) \deq \lim_{t \uparrow s} g(t)$.

For $s \in \Xi$, both $g(s-)$ and $g(s)$ are zeros of the function 
$x \mapsto \rho'(x) + \log(1 - x)$. For the sake of simplicity of
exposition, we shall assume below that there are never any zeros
of this function strictly between $g(s-)$ and $g(s)$:
\begin{equation} \label{eq:onlytwozeros}
  \{x \,:\,  s \rho'(x) + \log(1 - x) = 0 \}
  \bigcap (g(s-), g(s)) = \emptyset \,,
  \quad \text{for all } s \in \Xi \,.
\end{equation}
Also assume that $\Xi$ has no accumulation points. This is true, for
example, if $\sum_k k^2 \rho_k < \infty$.

Let $\{B_s, \ s \in \Xi\}$ denote a collection of
independent Bernoulli($1/2$) random variables, indexed by the
discontinuity set \eqref{eq:discontinuity_set}.
Define 
\begin{equation} \label{eq:ifl}
  \begin{split}
    \tilde{T}_t & \deq g(t-) +  B_t(g(t) - g(t-)) \,, \quad t \in \Xi \\
    \tilde{T}_t & \deq g(t), \ t \not\in \Xi \,.
  \end{split}
\end{equation}
In other words, at each point of discontinuity we choose the left
limit or the right limit of $g$ according to the flip of a fair coin. 
Finally, let
\begin{equation} \label{eq:z_tilde_definition}
  \tilde{Z}_t \deq t\rho(\tilde{T}_t) 
  - (1 - \tilde{T}_t) \log(1 - \tilde{T}_t) \,.
\end{equation}

For a sequence of stochastic processes $\{X^N\}_{N=1}^\infty$,
where $X^N = \{X^N_t\}_{t \geq 0}$, and a
stochastic process $X = \{X_t\}_{t \geq 0}$, we write $X^N
\convfdd X$ if the finite-dimensional distributions of $X^N$
converge to those of $X$.  For a sequence of random variables
(or vectors) $\{X^N\}$, we write $X^N \weakconv X$ to indicate that $X^N$
converges in distribution to $X$.
\begin{mainthm} \label{thm:intro_lambda_geq_one}
  Consider a Poisson hypergraph process such that $\rho_1 > 0$, and
  suppose \eqref{eq:onlytwozeros} holds. As $N \rightarrow \infty$, 
  \begin{equation} \label{eq:fluid_limit_with_patches}
    \{(\tilde{T}^N_t, \tilde{Z}^N_t)\}_{t \geq 0}
    \convfdd 
    \{(\tilde{T}_t, \tilde{Z}_t)\}_{t \geq 0} \,.
  \end{equation}
  Furthermore for any compact interval $I \subset [0,\infty) \setminus
  \Xi$, 
  \begin{equation} \label{eq:convinprob}
    \sup_{t \in I}\left|(\tilde{T}^N_t, \tilde{Z}^N_t) - 
    \left( g(t), t\rho(g(t)) - [1 - g(t)]\log(1 - g(t)) \right) \right|
    \rightarrow 0
  \end{equation}
  in probability as $N \rightarrow \infty$.
\end{mainthm}

We now turn to the case of patch-free hypergraph processes,
i.e.\ the regime where $\rho_1  = 0 < \rho_2$. 
%By virtue of \eqref{eq:rho1_plus_rho2}, 
In this case $g(s) = 0$ for all 
$s \in  [0, (2\rho_2)^{-1})$. There are three
possibilities for the behavior of $g$ at 
$(2\rho_2)^{-1}$, enumerated in Table \ref{tab:rho_1_equals_zero}.

\begin{table}[h]
  \begin{center}
    \begin{tabular}{cc}
      \textsc{Sub-case of $\rho_1 = 0 < \rho_2$} & \textsc{
	Behavior of $g$} \\
      \hline
      $3\rho_3 < \rho_2$ & \parbox{2.0in}{\vspace{0.1in} $g$ is
	continuous at $(2 \rho_2)^{-1}$, \\ and right derivative is
	finite \vspace{0.1in}} \\ 
      \hline
      $3\rho_3 = \rho_2$ & \parbox{2.0in}{\vspace{0.1in}
	$\rho_4,\rho_5,\ldots$ determine whether \\ 
	$g$ is continuous at $(2 \rho_2 )^{-1}$ \vspace{0.1in}} \\
      \hline
      $3 \rho_3 > \rho_2$ & \parbox{2.0in}{\vspace{0.1in}
	$g$ is discontinuous at $(2 \rho_2)^{-1}$}
    \end{tabular}

    \vspace{0.1in}

    \caption{The $0 = \rho_1 < \rho_2$ regime. \label{tab:rho_1_equals_zero}} 
  \end{center}
\end{table}

For simplicity, we focus on the case where $g$ has a single
discontinuity, located at $(2\rho_2)^{-1}$; i.e.\ 
$\Xi = \{(2\rho_2)^{-1}\}$. The general case follows the same pattern as
Theorem \ref{thm:intro_lambda_geq_one}, because
after the number of identifiable vertices has reached $O(N)$, the
subsequent evolution is much the same as the $\rho_1 > 0$ case. 

In the $\rho_1 = 0$ and $\rho_2 > 0$ regime, another  structure function
besides \eqref{eq:structure_function} comes into
play, namely the structure function $t_2$ of the graph which
results from discarding all hyperedges of cardinality more than two: 
\begin{equation} \label{eq:t2}
  t_2(x) \deq  \frac{-\log(1 - x)}{2 \rho_2x} \,,  
  \quad  x \in (0, 1) \,.
\end{equation}
Since $t_2$ is monotonic, the corresponding lower envelope 
$g_2$ defined as
\begin{equation} \label{eq:g2s}
  g_2(s) \deq \inf\{ x \in (0, 1) \,:\, t_2(x) > s \}\,,  
  \quad s \geq  0\,, 
\end{equation}
is continuous.  As before, $g_2(s) = 0$ for $0 \leq  s \leq
(2\rho_2)^{-1}$, and $g_2(s) \rightarrow 1$ as $s\rightarrow \infty$;
it describes the asymptotic proportion of vertices in the giant
component of a random graph where the ratio of edges to vertices is
$s\rho_2$.   

We will construct in Section \ref{sec:patch-free_process} an increasing
process $\{M_t\}$ so that the distribution of $M_t$ is 
\begin{equation} \label{eq:bp_population_0}
  \Pb{M_t = n} =
  \begin{cases} 
    e^{-2 \beta_2 n} \left( 2 t \rho_2 n \right)^{n - 1} / n! 
           & \text{if } n \in \N \,, \\
    \phi_t & \text{if } n = \infty \,,
  \end{cases} 
\end{equation}
where $\phi_t$ is the largest solution $x$ in $[0,1]$ of $2 t \rho_2 x + 
\log(1 - x) = 0$. [Notice that $\phi_t = 0$ for $2 t \rho_2 \leq 1$, and 
$0 < \phi_t < 1$ otherwise.]

Write $T^N_t$ for the number of vertices in the domain of $v_0$ in
$\Lambda_t$, and write $Z^N_t$ for the number of hyperedges
identifiable from $v_0$ in $\Lambda_t$. Set  
$\bar{T}^N_t \deq N^{-1}T^N_t$ and 
$\bar{Z}^N_t \deq N^{-1}Z^N_t$. 
Also, define
\begin{align*}
  \bar{T}_t & \deq g(t) \one_{\{M_t = \infty\}} \,; \\ 
  \bar{Z}_t & \deq \left\{ t\rho(g(t)) - [1 - g(t)]\log(1 - g(t))
              \right\} \one_{\{M_t = \infty\}} \,.
\end{align*}
\begin{mainthm} \label{thm:intro_no_patch_process}
  Consider a Poisson hypergraph process such that
  $\rho_1  = 0 < \rho_2$, and suppose $g$ has a
  single discontinuity located at $(2\rho_2)^{-1}$. Fix a distinguished
  vertex $v_0$. The number of vertices in the domain of $v_0$, and
  number of hyperedges identifiable from $v_0$, obey the following limits
  in distribution as $N \rightarrow \infty$: 
  \begin{equation} \label{eq:micro_no_patch_fluid_limit}
    \{(T^N_t, Z^N_t)\}_{t \geq  0} 
    \text{ converges weakly in $D\left([0, \infty), (\N \cup
        \{\infty\})^2 \right)$ to }  
      \{(M_t, M_t)\}_{t \geq  0} \,,
  \end{equation}
  where we adjoin $\infty$ to $\N$ as a
  compactifying point.   Also 
    \begin{equation} \label{eq:macro_no_patch_fluid_limit}
      \{(\bar{T}^N_t, \bar{Z}^N_t)\}_{t \geq  0}
      \convfdd \{ (\bar{T}_t,\bar{Z}_t) \}_{t \geq  0} \,.
    \end{equation}
\end{mainthm}

\begin{rmk}
  Observe the difference between the limit law
  $\{(\bar{T}_t,\bar{Z}_t)\}_{t \geq  0}$ in
  \eqref{eq:macro_no_patch_fluid_limit} and the limit law $\{
  (\tilde{T}_t, \tilde{Z}_t)\}_{t \geq  0}$ in
  \eqref{eq:fluid_limit_with_patches}:  $\tilde{T}_t$ conforms to the
  deterministic lower envelope $g(t)$, except at points in the finite
  discontinuity set, whereas $\bar{T}_t$ waits until the random time
  $\chi \deq \inf\{t \geq  0 \,:\, M_t = \infty\}$, with distribution
  function $g_2(t)$, before jumping from $0$ up to $g(t)$.  
\end{rmk}
\begin{rmk}
  See Remark \ref{rmk:no_conv_in_D} as to whether the convergence
  \eqref{eq:macro_no_patch_fluid_limit} extends to weak convergence in
  the Skorohod space $D([0, \infty),\R_+^2)$. 
\end{rmk}

The rest of this paper is organized as follows: Some definitions
concerning hypergraphs are given in Section \ref{sec:definitions}. We
establish that certain key processes are Markov in Section
\ref{sec:markov}.  The case of hypergraphs and hypergraph processes
with patches are treated in Section \ref{sec:review} and Section
\ref{sec:patches} respectively. Theorem \ref{thm:intro_lambda_geq_one}
is proved in Section \ref{sec:patches}.  Patch-free random hypergraphs
and hypergraph processes are treated in Section
\ref{sec:patch-free_static} and Section \ref{sec:patch-free_process},
respectively. Theorem 2 is proved in Section
\ref{sec:patch-free_process}. Finally, we mention future directions in
Section \ref{sec:future}. 

\section{Hypergraph definitions} \label{sec:definitions} 

Recall from the Introduction that the identifiable vertices are those vertices
removed by the hypergraph collapse algorithm described there, and the
identifiable hyperedges are those hyperedges consisting only of
identifiable vertices. 

Given a hypergraph $\Lambda$ and a subset $S \subset V$,
$\Lambda^S$ denotes the hypergraph after all vertices in $S$ are
deleted.  More precisely,
\begin{equation} \label{eq:LambdaS}
  \Lambda^S(A) \deq \sum_{B \supset A, \ B\setminus S = A} 
  \Lambda(B), \quad A \subset V\setminus S \,.
\end{equation}
We now more exactly specify the hypergraph collapse algorithm: select
if possible a vertex $v$ with 
$\Lambda(\{v\}) \geq 1$; replace $V$ by $V \setminus \{v\}$ and
$\Lambda$ by $\Lambda^{\{v\}}$; then repeat. When
the algorithm terminates, we obtain a set
$V^\star$ consisting of the identifiable
vertices, and a patch-free hypergraph $\Lambda^{V^\star}$ on
$V \setminus V^\star$.

Suppose $\Lambda$ is a patch-free hypergraph, and thus 
having no identifiable vertices. Given such
a hypergraph $\Lambda$ and a distinguished vertex $v_0$, we
say that $v$ is in the \emph{domain} of 
$v_0$ in  $\Lambda$ if $v$ is identifiable in the hypergraph 
$\Lambda + \one_{\{v_0\}}$ obtained by augmenting $\Lambda$
by the hyperedge $\{v_0\}$. A hyperedge is said to be 
\emph{identifiable} from  $v_0$ if it is identifiable in $\Lambda +
\one_{\{v_0\}}$. 

\begin{figure}
  \begin{center}
    \includegraphics{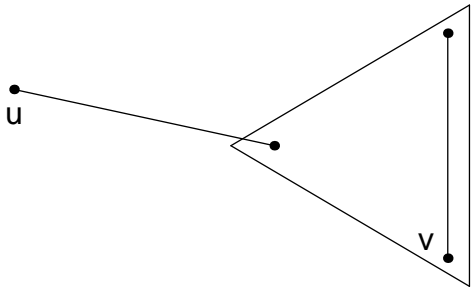}
    \caption{Adding a patch on $v$ makes $u$ identifiable, but not 
      vice versa. \label{fig:identifiable}}
  \end{center}
\end{figure}

\emph{Warning}: For a general patch-free hypergraph, it is possible for
vertex $u$ to be in the domain of $v$, while $v$ is not in the domain of $u$,
although this cannot happen in graphs; see Figure
\ref{fig:identifiable}.

\section{Poisson Hypergraph Processes: Markov Properties}
\label{sec:markov}

For $\{\beta_k\}_{k=1}^\infty$ a sequence of non-negative
numbers,
a \emph{Poisson($\beta$) random hypergraph} is a random hypergraph 
$\Lambda$ with vertex set $V = \{1,\ldots,N\}$ so
that for $A \subset V$, 
\begin{enumeratei}
  \item the random variable $\Lambda(A)$ has
    a Poisson distribution with mean $N \beta_{|A|}/\binom{N}{|A|}$, and
  \item  $\{ \Lambda(A) \,:\, A \subset V\}$ is
    a collection of independent random variables.
\end{enumeratei}

In what follows, $\{\rho_k\}_{k=1}^\infty$ will be a probability
distribution on the positive integers which has finite mean and
\begin{equation} \label{eq:rho1_plus_rho2}
  \rho_1 + \rho_2 > 0 \,.
\end{equation}
We now give an explicit construction of the hypergraph-valued
stochastic process described in the introduction.
Let $K_1, K_2, \ldots$ be a sequence of independent random variables
in $\{1,2,3,\ldots\}$ with common distribution
$\Pb{K_n = k} = \rho_k$, for all $n,k \in \N$.  Denote by
$A_1, A_2, \ldots$ a sequence of independent random subsets of $V$,
such that $A_n$ is chosen uniformly at random from the subsets of $V$ of
size $K_n$ whenever $K_n \leq N$; the set $A_n$ is not defined when
$K_n > N$. Let $\{E_t\}_{t \geq 0}$ be a Poisson process, run at rate
$N$, having arrival times $\tau_1, \tau_2, \ldots$. Define a stochastic
process $\{\Lambda_t\}_{t \geq 0}$ with values in the set of
hypergraphs with vertex set $V$ by 
\begin{equation*}
  \Lambda_t(A) \deq \sum_{n \,:\, \tau_n \leq t}\one_{\{A = A_n\}} \,.
\end{equation*}
Interpret $\Lambda_t(A)$ as the number of
occurrences of hyperedge $A$ by time $t$. 
In summary, for each $A \subset V$, 
\begin{equation} \label{eq:pprate}
  \{\Lambda_t(A)\} \text{ is a Poisson process of rate }
  N \frac{\rho_{|A|} }{ \binom{N}{|A|} } \,,
\end{equation}
and all these Poisson processes are independent. We call
$\{\Lambda_t\}_{t \geq 0}$ a \emph{Poisson($\rho$) hypergraph process},
where $\rho$ denotes the generating function 
\begin{equation}
  \rho(x) \deq \sum_{k \geq 1} \rho_k x^k \,.
\end{equation}
The finite mean assumption is equivalent to $\rho'(1) < \infty$.  
For fixed $t \geq 0$, $\Lambda_t$ is a Poisson($t \rho$) random
hypergraph. 

Whereas the hypergraph literature has tended to concentrate on
the ``$k$-uniform'' case (i.e.\ $\rho_k = 1$ for some $k$), we
find the superposition of $k$-uniform random hypergraphs for
various different values of $k$ can be handled without special effort,
and leads to asymptotic properties absent from the $k$-uniform
case.  Moreover the Poisson structure simplifies our arguments, for
example by allowing some summary statistics of $\{\Lambda_t\}_{t \geq
  0}$ to be Markov processes in their own right: see Proposition 
\ref{prop:markov_property}.  Poissonization is, of course,
a well-established procedure -- see \ocite{A:PCH}.

Previous literature has also concentrated on the case $\Lambda\le1$. 
We now sketch a way to deduce from our results for a Poisson($\beta$) random 
hypergraph $\Lambda$ some corresponding results for $\Lambda \wedge 1$.
We note moreover that if $\rho_k = 1$ for some $k$ then $\Lambda \wedge 1$
is exactly a $k$-uniform hypergraph.
The set of identifiable vertices is the same for $\Lambda$ and
$\Lambda \wedge 1$ but $\Lambda$ may have additional identifiable hyperedges.
First consider patches. Throwing a Poisson($N\beta$) number of balls
(i.e.\ patches) uniformly at random into $N$ urns yields a
Binomial($N$,$1 - e^{-\beta_1}$) number of occupied urns (i.e.\
vertices covered by at least one patch).  Hence the number of patches
in $\Lambda$, less the number in $\Lambda \wedge 1$, divided by $N$,
has limit in probability $\beta_1 + e^{-\beta_1} - 1$.  On the other hand,
the expected number of subsets of size at least 2 receiving at least 2
hyperedges is bounded uniformly in $N$.
Hence, after rescaling by $N^{-1}$, only the extra patches in $\Lambda$
can contribute in the limit and of course all of these do so.

\begin{prop} \label{prop:markov_property}
   Let $T_t$ and $Z_t$ denote the numbers of identifiable vertices and
   identifiable hyperedges for $\Lambda_t$. Both $\{T_t\}_{t \geq 0}$ and
   $\{(T_t,Z_t)\}_{t \geq 0}$ are Markov processes. The
   number of non-identifiable hyperedges in $\Lambda_t$, given
   that $T_t = m$, is conditionally Poisson,
   with mean 
   \begin{equation} \label{eq:non_identifiable_exactly}
     N t \left[ 1 - \sum_{k \geq 1} \rho_k 
       \frac{\binom{m}{k} + (N-m) \binom{m}{k-1}}
	    {\binom{N}{k}} \right] \,.
   \end{equation}
   When $m - N\gamma = o(N)$,  for $\gamma \in [0,1]$, this reduces as $N
   \rightarrow \infty$ to 
   \begin{equation} \label{eq:nonidentifiable}
      N t \left[1 - \rho(\gamma) - (1-\gamma)\rho'(\gamma) \right]
      + o(N) \,.
   \end{equation}
\end{prop}
\begin{rmk}
  Because the total number of hyperedges in $\Lambda_t$ is
  Poisson($Nt$), %\eqref{eq:nonidentifiable} 
  Proposition \ref{prop:markov_property} reduces the study of limits
  of identifiable hyperedges to study of limits of identifiable
  vertices. In particular, if $N^{-1} T_t$ converges in distribution
  as $N \rightarrow \infty$ to a random variable $\tilde{T_t}$, then
  necessarily  
  \begin{equation} \label{eq:zt}
    N^{-1} Z_t \weakconv t\left[ \rho( \tilde{T}_t )
      + (1 - \tilde{T}_t) \rho'(\tilde{T}_t) \right] \,.
  \end{equation}
\end{rmk}
\begin{rmk}  
  It is easy to identify the generator of $\{T_t\}_{t \geq 0}$,
  rescale by division by $N$, and take a limit on any compact
  interval $I \subset \R_+ \setminus \Xi$ (see
  \eqref{eq:discontinuity_set}); however this approach did not lead to
  a proof of Theorem \ref{thm:intro_lambda_geq_one},  because of the
  difficulty of passing through discontinuous phase transitions.   
\end{rmk}
To prepare for the proof, some measure-theoretic
apparatus is needed. Let $(\Omega,\F,\P)$ be the
probability space on which the process $\{\Lambda_t\}_{t \geq 0}$
is defined.
For any set $S \subset V$, and any $t \geq 0$,
define the $\sigma$-field $\F_t^S \subset \F$ as
\begin{equation*}
  \F_t^S \deq \bigvee_{0 \leq s \leq t} 
  \sigma\{ \Lambda_s(A) \; : \; |A \setminus S| \leq 1 \} \,.
\end{equation*}
Let $V_t^\star$ denote the set of vertices identifiable at time $t$. By
construction, the event $\{ V_t^\star = S \}$ occurs if and only if,
%among supersets of the union of all vertices covered by patches, 
among all sets containing all vertices covered by patches,
$S$ is
the minimal subset of $V$ for which $\Lambda_t(A) = 0$ 
whenever $|A \setminus S| = 1$. Thus $\{V^\star_t = S\} \in \F_t^S$.

When we consider $V^\star_t$ as a ``stopping set'' for a 
set-indexed process, it becomes natural to define another 
$\sigma$-field:
\begin{equation*}
  \F_{V^\star_t} \deq \left\{ B \in \F \,:\,
  B \cap \{V^\star_t = S\} \in \F_t^S 
  \text{ for all } S \subset V \right\} \,.
\end{equation*}
$T_s$ and $Z_s$ are $\F_{V^\star_t}$-measurable, for all
$0 \leq s \leq t$.  We may describe $\F_{V^\star_t}$ informally
as the knowledge we have about $\{\Lambda_s\}_{0 \leq s \leq t}$
after performing hypergraph collapse at each time $s \in [0,t]$.

\begin{lem} \label{lem:stopset}
\hspace{0.1in}

\begin{enumeratei}
\item \label{it:stopset1}
  Fix any $t > 0$. Pick any collection of non-negative integers 
  $\{k_A \,:\, A \subset V \}$,
  and set
  \begin{equation*}
    p(S) \deq \P\!\left( \bigcap_{ A \,:\, |A \setminus S| > 1}
    \left\{ \Lambda_t(A) = k_A \right\} \right) \,.
  \end{equation*}
  Then
  \begin{equation*}
    \P\!\left( \bigcap_{A \,:\, |A \setminus V^\star_t| > 1}
    \left\{ \Lambda_t(A) = k_a \right\} \Given
    \F_{V^\star_t} \right) = p(V^\star_t) \,.
  \end{equation*}
\item \label{it:stopset2}
  Fix any $t > 0$. The conditional distribution of the random hypergraph
  $\Lambda^S_t$ (in the notation of \eqref{eq:LambdaS}),
  given $\F_{V^\star_t}$, on the event $\{V^\star_t = S\}$, where
  $|S| = m$, is that of a Poisson($\beta$) random hypergraph on $N - m$
  vertices with parameters 
  \begin{align}
    \beta_1 & \deq 0 \nonumber \\
    \beta_j & \deq \frac{t}{1 - m/N}
    \binom{N-m}{j} \sum_{i \geq 0} \rho_{i+j}
    \frac{\binom{m}{i}}{\binom{N}{j+i}}  \,, \quad j \geq 2
    \,. \label{eq:betaj} 
  \end{align}
\end{enumeratei}
\end{lem}
For a random variable $X$, we write $X \sim \text{Poisson}(\mu)$
to indicate that the distribution of $X$ is Poisson with
expectation $\mu$.  Also we will write $X \sim \text{Binomial}(n,p)$
to indicate that $X$ is a Binomial random variable with
parameters $n$ and $p$.
\begin{proof}[Proof of \eqref{it:stopset1}]
Certainly $p(V^\star_t)$ is $\F_{V^\star_t}$-measurable. It remains to
show that, for any $B \in \F_{V^\star_t}$, 
\begin{equation*}
  \int_B p(V^\star_t) d\P = \P\!\left( B \cap \bigcap_{A \,:\, 
    |A \setminus V^\star_t| > 1} \{ \Lambda_t(A) = k_a \} \right) \,.
\end{equation*}
Split the event on the right into disjoint events by intersecting with
$\{V^\star_t = S\}$ for each $S \subset V$. For each $S$,
$B \cap \{V^\star_t = S \}$ lies in $\F^S_t$, and therefore
is independent of $\{\Lambda_t(A) = k_a\}$ for every $A$ such that
$|A \setminus S| > 1$, by construction of a Poisson hypergraph
process.  The right side becomes
\begin{equation*}
  \sum_{S \subset V} p(S) \P\!\left(B \cap \{V^\star_t = S \}\right)
\end{equation*}
which is equal to the left side; \eqref{it:stopset1} follows.
\end{proof}

\begin{proof}[Proof of \eqref{it:stopset2}]
Suppose $S \subset V$ and $A \subset V \setminus S$ with
$|A| = j \geq 2$.  For any $C \subset S$ with $|C|=i$,
%Condition on $\{V^\star_t = S\}$, and fix $A \subset V \setminus S$
%with $|A| = j \geq 2$.  For any $C \subset S$ with $|C|=i$,
\eqref{eq:pprate} implies that
\begin{equation*}
  \Lambda_t(A \cup C) \sim
	 {\rm Poisson}\left( t \rho_{j+i} N / \binom{N}{j+i} \right) \,.
\end{equation*}
The result of part \eqref{it:stopset1} implies that
the random variables $\Lambda_t(A \cup C)$
are conditionally independent for different choices of
$C$, given $\{V^\star_t = S \} \cap \F_{V^\star_t}$.

If $|S| = m$, there are $\binom{m}{i}$ choices of $C$, and
following the notation of \eqref{eq:LambdaS},
\[
\Lambda^S_t(A) = \sum_{C \subset S} \Lambda_t(A \cup C)
   \sim {\rm Poisson}\left( t N \sum_{i \geq 0} \rho_{j+i}
   \binom{m}{i}/\binom{N}{j+i} \right) \,.
\]
In a Poisson($\beta$) random hypergraph on $(N-m)$ vertices, the
number of occurrences of $A$, where $|A|=j$, is Poisson with parameter
\[
(N-m) \beta_j/\binom{N-m}{k} \,.
\]
On comparison with the previous line, this verifies the formula 
\eqref{eq:betaj} for $\beta_j$, when $j \geq 2$.  Clearly there
are no $1$-hyperedges in $\Lambda^S_t$ when
$\{V^\star_t = S \}$, by definition of identifiability.  Hence
\eqref{it:stopset2} is established.
\end{proof}

\begin{proof}[Proof of Proposition \ref{prop:markov_property}]
Fix any $t > 0$. Suppose that $T_t = m$. The first jump in 
the process $\{(T_s,Z_s)\}_{s \geq t}$ can occur only when a new hyperedge
arrives, and the arrival time is independent of the past. The law of
the jump depends only on two things: the set $A$ of vertices in the new
hyperedge (which is independent of the past), and on the hypergraph
$\Lambda^S_t$, where $S \deq V^\star_t$.
Lemma \ref{lem:stopset}\eqref{it:stopset2}
establishes that the law of $\Lambda^S_t$,
conditional on $\F_{V^\star_t}$ is  fully determined
by $m$, $t$, and the parameters $\{\rho_i\}_{i \geq 1}$;
in particular it is conditionally independent of 
$\{(T_s,Z_s)\}_{0 \leq s \leq t}$
given that $\{T_t = m\}$. Hence the Markovian
property of $\{T_t\}_{t \geq 0}$ and
$\{(T_t,Z_t)\}_{t \geq 0}$ is established.

It follows from Lemma \ref{lem:stopset} that the total number of non-identifiable
hyperedges in $\Lambda_t$, given that $\{T_t = m\}$, is conditionally
Poisson, with mean $(N - m) \sum \beta_j$,
for $\beta_j$ as in \eqref{eq:betaj}.
Write  $k \deq i+j$, and switch the order
of summation, to obtain 
\[
\left( 1 - \frac{m}{N} \right) \sum \beta_j
   = t \sum_{k \geq 2} \rho_k \sum_{j=2}^k
     \binom{N-m}{j}\binom{m}{k-j} / \binom{N}{k} \,.
\]
On considering the Hypergeometric($(N,N-m,k)$) distribution,
we see that the inner sum is
\[
1 - \left[ \binom{m}{k} + (N-m) \binom{m}{k-1} \right]
   / \binom{N}{k} \,.
\]
The last expression is zero when $k=1$, so $(N-m)\sum \beta_j$ takes
the form \eqref{eq:non_identifiable_exactly}.  When
$m - N \gamma = O(N)$, the last expression converges,
as $N \rightarrow \infty$, to $1 - \gamma^k -
k\gamma^{k-1}(1-\gamma)$, and is bounded between $0$ and $1$.
The Bounded Convergence Theorem yields \eqref{eq:nonidentifiable}.
\end{proof}

\section{Identifiability In Random Hypergraphs With Patches}
\label{sec:review}
In this section we review some material from \ocite{DN:HG}. 

Fix $t > 0$, and set $\Lambda \deq \Lambda_t$, $\beta_k 
\deq t\rho_k$. In this case, $\Lambda$ is a Poisson($\beta$) random
hypergraph. Suppose we perform hypergraph collapse, described
above, in the following special way: at each step the next vertex $v$ to
be deleted is selected with a probability proportional to the number
of patches on $v$. This is called \emph{randomized collapse}. 
The \emph{debris} of a hypergraph is the number of hyperedges
equal to the empty set.
Set
$\Lambda_0 \deq \Lambda$, and let 
$\{\Lambda_n\}_{n \in \N}$ denote the sequence of hypergraphs
obtained. Set $Y_n$ and $Z_n$ to be the amount of patches and debris,
respectively, in $\Lambda_n$; formally 
\[
Y_n \deq \sum_{v \in V} \Lambda_n(\{v\})\,, \quad \text{and} \quad 
Z_n \deq \Lambda_n(\emptyset) \,.
\]
The key observation in \ocite{DN:HG} is that 
$\{(Y_n,Z_n)\}_{n \in \N}$ is a Markov chain (but not the same one as in
Proposition \ref{prop:markov_property}, for here $t$
is fixed!), which stops at
\begin{equation} \label{eq:T_definition} 
  T \deq \inf\{n \,:\, Y_n=0\} \,.
\end{equation}
Moreover, conditional on $\{Y_n = m, Z_n = k\}$, 
\begin{equation} \begin{split} \label{eq:y_and_z}
    Z_{n + 1} & = k + 1 + W_{n + 1} \,, \\
    Y_{n + 1} & = m - 1 - W_{n + 1} + U_{n + 1} \,.
\end{split} \end{equation}
Here $W_{n+1}$ and $U_{n+1}$ are independent, with
\begin{equation} \begin{split} \label{eq:u_and_w}
    W_{n + 1} & \sim {\rm Binomial}\left( m - 1, \frac{1}{N-n} \right) \\
    U_{n + 1} & \sim {\rm Poisson}\left( (N-n-1)t\lambda_2(N,n) \right) 
\end{split} \end{equation}
where
\begin{equation} \label{eq:lambda2}
  \lambda_2(N,n)  \deq N \sum_{i=0}^n \rho_{2 + i}
  \binom{n}{i} / \binom{N}{i+2} \,.
\end{equation}

By construction, $T = |V^\star|$, the number of
identifiable vertices, and $Z \deq Z_T \deq \Lambda_T(\emptyset)$
is the number of identifiable hyperedges. For comparison, note that, by
%\eqref{eq:nonreducible}, 
Proposition \ref{prop:markov_property}
the number of
non-identifiable hyperedges in $\Lambda$, given that
$T = N\gamma$, is conditionally Poisson, with mean 
\begin{equation} \label{eq:poisson_mean_asymp}
N(t - \beta(\gamma) - 
   (1 - \gamma)\beta'(\gamma)) + o(N) \,.
\end{equation}

%By passing to the limit as $N \rightarrow \infty$ for the Markov
%chain $\{Y_n\}_{n \in \N}$ after rescaling, 
We obtained a limit theorem for
$\tilde{T}^N \deq N^{-1}T$ and $\tilde{Z}^N \deq N^{-1}Z$,
where $Z$ is the number of identifiable hyperedges.  We state the result
in a simple case. Set  
\[
\beta(x) \deq \sum_k \beta_k x^k \,, \quad
x \in [0, 1] \,.
\]
Assume that $\beta_1 > 0$ and that the derivative $\beta'(1)
<\infty$. 
Then
\begin{equation} \label{eq:18}
  \{x \in [0, 1) \,:\, \beta'(x) + \log(1 - x) < 0 \} 
\end{equation}
is non-empty, and its infimum is $g(t)$, as defined in
\eqref{eq:gsdefinition}. By our assumption
\eqref{eq:onlytwozeros}, there is at most
one $x \in [0,g(t))$ such that $\beta'(x) + \log(1 - x) = 0$,
namely $g(t-)$; this is different to $g(t)$ only if
$t \in \Xi$, the set of discontinuity points of the lower envelope
$g$.  

Let $\tilde{T}$ be a random variable taking values $g(t)$ and 
$g(t-)$, each with probability $1/2$. As a special case of 
of \ocite{DN:HG}*{Theorem 2.2} we know: 

\begin{thm} \label{thm:fixedtime}
  The following limit in distribution holds as $N \rightarrow \infty$:
  \begin{equation} \label{eq:tandzlimits}
    \left( \tilde{T}^N,\tilde{Z}^N \right) \weakconv
    \left( \tilde{T}, \beta(\tilde{T}) - (1 - \tilde{T})
    \log(1 - \tilde{T}) \right) \,.
  \end{equation}
\end{thm}

\begin{rmk}
  \ocite{GN:EE} have shown that the limit for the rescaled number of
  identifiable hyperedges can be decomposed as follows: $(1 -
  \tilde{T})\log(1 - \tilde{T})$ counts the \emph{essential
  hyperedges}, i.e.\ those whose absence would have reduced the set of
  identifiable vertices, and $\beta(\tilde{T})$ counts the remainder. 
\end{rmk}

\begin{rmk}
  Suppose in particular that $\Lambda \deq \Lambda_t$
  and $\beta(x) \deq t\rho(x)$ for some $t \in \Xi$, the
  discontinuity set of $g$. Then \eqref{eq:tandzlimits} implies that
  the proportion of identifiable vertices has a limit in distribution
  which is random, taking the values $g(t)$ and $g(t-)$ each with
  probability $1/2$.   
\end{rmk}

\begin{rmk}
  It suffices to derive the limit for $\tilde{T}^N$,
  since the limit for $\tilde{Z}^N$ follows from Proposition
  \ref{prop:markov_property}.  To check this, recall that, by
  \eqref{eq:zt}, if $\tilde{T}^N$ converges to $g(t)$, then the number
  of identifiable hyperedges, divided by $N$, converges to  
  \begin{equation} \label{eq:formula1}
    t \left\{ \rho( g(t) ) + [1 - g(t)]\rho'( g(t) ) \right\} \,.
  \end{equation}
  However by definition of $g(t)$, 
  $t\rho'(g(t)) = -\log(1 - g(t))$, so we have recovered the 
  formula $\beta(\tilde{T}) - (1 - \tilde{T})\log(1 - \tilde{T})$.
\end{rmk}

\section{Identifiability In Hypergraph Processes With Patches}
\label{sec:patches}

In this section we move from the static random hypergraph model of Theorem
\ref{thm:fixedtime} to the Poisson($\rho$) hypergraph process
$\{\Lambda_t\}_{t \geq 0}$, providing here a proof of
Theorem \ref{thm:intro_lambda_geq_one}.

Extending the notation of the previous section, let $\tilde{T}^N_t$
and $\tilde{Z}^N_t$ denote the rescaled numbers of identifiable
vertices and hyperedges for $\Lambda_t$, respectively, as
defined in \eqref{eq:T_tilde_N_and_Z_tilde_N_defn}. Note that $t
\mapsto \tilde{T}^N_t$ and $t \mapsto \tilde{Z}^N_t$ are increasing,
right-continuous, stochastic processes. It follows from Proposition
\ref{prop:markov_property} that $\{(\tilde{T}^N_t, \tilde{Z}^N_t)\}_{t
  \geq 0}$ is a Markov process.  

\begin{proof}[Proof of Theorem \ref{thm:intro_lambda_geq_one}]
Fix $0 \leq t_1 < \ldots < t_r$. We have to show the 
convergence in distribution
\begin{equation} \label{eq:tz_limit_law}
\{ (\tilde{T}^N_{t_i}, \tilde{Z}^N_{t_i}) \}_{i = 1, \ldots, r} 
\weakconv \{ (\tilde{T}_{t_i}, \tilde{Z}_{t_i}) \}_{i = 1,
  \ldots, r}  \,.
\end{equation}
It suffices to do so when at least one of $\{t_i, t_{i + 1}\}$ is not a
discontinuity point, for every $i \in \{1,\ldots,r-1\}$. 
Proposition \ref{prop:markov_property} showed that
$\{(\tilde{T}^N_t,\tilde{Z}^N_t)\}_{t \geq 0}$ is Markov,
and for any Markov process $\{Y_t\}_{t \geq 0}$ the conditional
law of $Y_{t_r}$ given $(Y_{t_1},\ldots,Y_{t_{r-1}})$ is the same
as the conditional law given $Y_{t_{r-1}}$.  Hence it suffices
to consider the case $r=2$ such that $t_1 \not\in \Xi$ or
$t_2 \not\in \Xi$, and these possibilities are both subsumed in the
case $r=3$ with $t_1,t_3 \not\in \Xi$.
Then only the marginal limit at time
$t_2$, as given in Theorem \ref{thm:fixedtime} is random,
so Theorem \ref{thm:fixedtime} implies the
full convergence in distribution. 

The second assertion follows from the first since all processes are
increasing, and the limit is deterministic and continuous on
$I$.
\end{proof}

\begin{rmk} \label{rmk:no_conv_in_D} 
  The rescaled number of essential hyperedges, as studied by 
  \ocite{GN:EE}, has a limit  
  $\{-(1 - \tilde{T}_t)\log(1 - \tilde{T}_t) \}_{t \geq 0}$ in the
  same sense as \eqref{eq:fluid_limit_with_patches} and
  \eqref{eq:convinprob}. 
\end{rmk}

\begin{rmk}
  One may ask whether the convergence
  \eqref{eq:fluid_limit_with_patches} extends to weak convergence in
  the Skorohod space $D([0,\infty),\R_+^2)$. Since $t \mapsto
  \tilde{T}^N_t$ and $t \mapsto \tilde{Z}^N_t$ are non-decreasing, the
  necessary and sufficient condition of \ocite{JS:LT1}*{p.\ 306}
  may be applied, which would require that the sum of squared jumps of
  $\{\tilde{T}^N_t\}$ converges in law to the sum of squared jumps of
  $\{\tilde{T}_t\}$, and similarly for
  $\{\tilde{Z}^N_t\}$. Unfortunately the techniques presented in this
  paper do not seem to be able to confirm this; indeed, it seems
  plausible that, for arbitrarily large $N$, and for $t \in \Xi$,
  there is a probability bounded away from zero that $\tilde{T}^N_s$
  makes more than one jump in going from $\approx g(t-)$ to $\approx
  g(t)$ at time $s \approx t$, and this would contradict the condition
  stated.  
\end{rmk}

\begin{rmk}
  If \eqref{eq:onlytwozeros} is false, one can reformulate the process 
  \eqref{eq:ifl}, by consulting \ocite{DN:HG}*{Theorem 2.2} and prove
  a corresponding version of Theorem \ref{thm:intro_lambda_geq_one}.
\end{rmk}

\section{Domain Of A Vertex In A Hypergraph Without Patches}
\label{sec:patch-free_static}

We revert to the fixed-time setting of Section
\ref{sec:review}. Suppose $\Lambda$ is a Poisson($\beta$) random 
hypergraph, such that 
\[
\beta_0 = \beta_1 = 0 < \beta_2\,, \quad
   \beta(x) \deq  \sum_{k \geq 2}\beta_k x^k\,, \ 
    x \in [0, 1]\,. 
\]

Fix a vertex $v_0$. Write $T^N$ for the number of vertices in the domain
of $v_0$, and write $Z^N$ for the number of hyperedges identifiable from
$v_0$. Set $\bar{T}^N \deq N^{-1}T^N$ and $\bar{Z}^N \deq N^{-1}Z^N$. Both
the \emph{microscopic variables} $(T^N, Z^N)$, and the 
\emph{macroscopic variables}
$(\bar{T}^N, \bar{Z}^N)$ have non-trivial limits as $N \rightarrow
\infty$, which we now describe. The coefficient $\beta_2$ plays a distinguished
role. 

\begin{lem}
  Let $\{\xi_n\}_{n \in \N}$  be a random walk on the integers, started
  at $\xi_0 = 1$, whose increments are of the form 
  $\xi_n  - \xi_{n - 1} = -1 + {\rm Poisson}(2\beta_2)$. Let $\phi$ be
  the largest root in $[0,1]$ of $2 \beta_2 x + 
  \log(1 - x) = 0$, so $\phi = 0$ for $2 \beta_2 \leq 1$, and 
  $0 < \phi < 1$ otherwise. Then the first passage time to $0$,
  \begin{equation} \label{eq:M_definition}
    M \deq \inf\{n \geq 0 \,:\, \xi_n = 0 \} \,,
  \end{equation}
  has the following distribution:
  \begin{equation}
    \begin{split} \label{eq:bp_population}
      \Pb{M = n}      & = e^{-2 \beta_2 n} \left(
                        2 \beta_2 n \right)^{n - 1} /
			n!\,, \quad n \in \N \,; \\
      \Pb{M = \infty} & = \phi \,,
    \end{split}
  \end{equation}
\end{lem}
\begin{rmk*} 
  $M$ is distributed as the total number of individuals in a
  branching process with one ancestor, and Poisson($2\beta_2$) offspring
  distribution. This distribution describes the sizes of small
  components in an Erd\H{o}s-R\'enyi random graph; see
  \ocite{B:RG}. 
\end{rmk*}

\begin{proof}
The fact that $\Pb{M = \infty} = \phi$ is an elementary fact from the
theory of branching processes. The formula for $\Pb{M=n}$ is a special
case of a formula of \ocite{D:TP}, which is proved in detail on p.\ 300
of \ocite{D:BP}.
\end{proof}

Assume that $\beta'(1) < \infty$. Then the set
\eqref{eq:18} is non-empty, and
its infimum is $g \deq g(t)$, as defined in
\eqref{eq:gsdefinition}.  Assume further that
$\beta'(x) + \log(1 - x) > 0$ for all $x \in (0,g)$. If either of these
assumptions fail, then the techniques of \ocite{DN:HG}, 
combined with some arguments given below, still establish the desired
asymptotics. We omit the details. 

Set 
\begin{equation} \begin{split}
    \bar{T} & \deq g \one_{\{M = \infty\}} \,; \\
    \bar{Z} & \deq \left[ \beta(g) - (1 - g)\log(1 - g) \right]
              \one_{\{M = \infty\}}\,.
\end{split} \end{equation}

\begin{thm} \label{thm:identif_sans_patches} 
  Consider a Poisson random hypergraph without patches, and fix a distinguished
  vertex $v_0$. The number of vertices in the domain of $v_0$, and number
  of hyperedges identifiable from $v_0$, obey the following limits in
  distribution as $N \rightarrow \infty$: 
  \begin{equation}
    (T^N, Z^N) \weakconv (M, M)\,; \quad 
    (\bar{T}^N, \bar{Z}^N) \weakconv (\bar{T}, \bar{Z})\,.
  \end{equation}
  Here $M$ is considered as a random variable taking values in the
  one-point compactification $\N \cup \{\infty\}$ of $\N$.
\end{thm}

\begin{proof}
\emph{Step I}. Set $\Lambda_0 \deq \Lambda + \one_{\{v_0\}}$, and let 
$\{\Lambda_n\}_{n \in \N}$ be a sequence of hypergraphs obtained
by randomized collapse. Denote by $Y^N_n$ and $Z^N_n$ the numbers of
patches and debris, respectively, in $\Lambda_n$. Then  
\[
T^N \deq \inf\{n \geq 0 \,:\, Y^N_n = 0\}\,;
\quad Z^N \deq Z^N_{T^N} \,. 
\]

We know that $\{(Y^N_n, Z^N_n)\}_{n \geq 0}$ is a Markov chain,
starting from $(1,0)$: the increments, conditional on 
$Y^N_n = m \geq 1$ and $Z^N_n = k$, are as given in
\eqref{eq:y_and_z} and \eqref{eq:u_and_w}.
 
For fixed $n \geq 0$ and $m \geq 1$, the random variable $W_{n+1}$
defined in \eqref{eq:u_and_w} converges to $0$ in
distribution as $N \rightarrow \infty$. Also 
\begin{equation} \label{eq:lambda2_limit}
(N - n - 1) \lambda_2(N, n) \rightarrow 2  \rho_2 \,.
\end{equation}
so the random variable $U_{n+1}$ defined in 
\eqref{eq:u_and_w}
converges to Poisson($2\beta_2$) in distribution as 
$N \rightarrow \infty$. Hence, for all $n \geq 0$, 
\[
\{(Y^N_j, Z^N_j)\}_{0 \leq j \leq n}
        \weakconv \{(\xi_j, j)\}_{0 \leq j \leq n} 
\]
which implies $(T^N, Z^N) \weakconv (M, M)$ as $N \rightarrow \infty$.  
If $2\beta_2 \leq 1$, then $\Pb{M = \infty} = 0$, 
so the proof is complete. It only remains to prove the second
convergence assertion in the case where $2\beta_2 > 1$, and 
$0 < \phi < 1$.

\medskip 
\emph{Step II}. Introduce an auxiliary time variable $t$, and let
$\{\nu_t\}_{t \geq 0}$ be a Poisson process of rate $N$. Set  
\begin{align*}
  \bar{Y}^N_t   & \deq N^{-1}Y^N_{\nu_t} \,, \\ 
  \bar{Z}^N_t   & \deq N^{-1}Z^N_{\nu_t} \,, \\
  \bar{\nu}^N_t & \deq N^{-1}\nu_t \,,       \\
  \tau^N        & \deq \inf\{t \geq 0 \,:\,
                  \bar{Y}^N_t = 0 \} \,. 
\end{align*}
With reference to \ocite{DN:HG}, set
\begin{align*}
  y(t) & \deq (1 - t)(\beta'(t) + \log(1 - t)) \,;\\
  z(t) & \deq \beta(t) - (1-t)\log(1 - t) \,.
\end{align*}
By Theorem 6.1 and Remark 6.2 of \ocite{DN:HG},
for all $\delta > 0$,
\begin{equation} \label{eq:old_fluid_limit_result}
  \limsup_{N \rightarrow \infty} \frac{1}{N}
  \log\left( \Pb{ \sup_{t \leq \tau^N}
    \left\| (\bar{\nu}^N_t,\bar{Y}^N_t, \bar{Z}^N_t) - 
    (t, y(t), z(t)) \right\| > \delta} \right) < 0 \,.
\end{equation}
Observe that $\bar{\nu}^N_{\tau^N} = \bar{T}^N$, which will have the
same limit in probability as does $\tau^N$. We will show that, for
all $\theta \in (\log(1-\phi)],0)$, there exists $\delta > 0$
and $N_0$ such that 
\begin{equation} \label{eq:unlikely_to_die_before_delta}
  \Pb{\bar{T}^N \leq \delta} \leq e^{\theta} \,,
    \quad \text{for all } N \geq N_0 \,.
\end{equation}
By \eqref{eq:bp_population} and the fact
that $T^N \weakconv M$, we know that, for all
$\delta > 0$ and all $\phi' > \phi$: 
\begin{equation*}
  \Pb{\bar{T}^N \leq \delta} \geq 1 - \phi'
\end{equation*}
for all sufficiently large $N$. Also from 
\eqref{eq:old_fluid_limit_result} we obtain, for all 
$\delta > 0$,
\begin{equation*}
  \Pb{\bar{T}^N \in (\delta, g - \delta) \cup (g + \delta, \infty)}
  \rightarrow 0 
\end{equation*}
as $N \rightarrow \infty$. Hence the claim that 
$(\bar{T}^N, \bar{Z}^N) \weakconv (\bar{T}, \bar{Z})$ will follow
as soon as we have proved \eqref{eq:unlikely_to_die_before_delta};
then \eqref{eq:old_fluid_limit_result} will strengthen this
to show $(\bar{T}^N,\bar{Z}^N) \weakconv (\bar{T},\bar{Z})$.

\medskip
\emph{Step III}. The remainder of the proof is to establish
\eqref{eq:unlikely_to_die_before_delta}. Given $Y^N_n = m \geq 1$, set 
\begin{align*}
  \Phi^N(m, n) & \deq 
  \E \exp\left\{\theta(-1 - W_{n + 1} + U_{n + 1}) \right\}  \\
  & = \exp\left\{-\theta + F\left(m - 1, \frac{1}{N - n},
  -\theta\right) + G((N - n - 1)\lambda_2(N, n), \theta) \right\} \,.
\end{align*}
where 
\begin{equation*}
  F(k, p, \theta) \deq k \log\left(1 - p + p e^{\theta} \right) \,; 
  \quad G(\mu, \theta) \deq \mu(e^{\theta} - 1) \,.
\end{equation*} 
Lemma 6.1 of \ocite{DN:HG} implies that
\begin{equation*}
  \sup_{n \leq N/2} \left|(N - n - 1)\lambda_2(N, n) - 
  \left(1 - \frac{n}{N} \right) \beta''(n/N) \right|
  \rightarrow 0 \,,
\end{equation*}
as $N \rightarrow \infty$. Since $\theta > \log(1 - \phi)$, there is 
$\bar{\phi} < \phi$ such that $\theta > \bar{\theta}
\deq \log(1 - \bar{\phi})$; by construction of $\phi$,  
$2\beta_2\bar{\phi} + \log(1 - \bar{\phi}) > 0$, so 
$2 \beta_2 (1 - e^{\bar{\theta}}) + \bar{\theta} > 0$; 
in other words,
\begin{equation*}
  \exp\{ - \bar{\theta} + G(2\beta_2, \bar{\theta}) \} < 1 \,.
\end{equation*}
We can therefore find $\delta > 0$ and $N_0$ such that
\begin{equation} \label{eq:phi_bound}
  \Phi^N(m, n) \leq 1\,, \quad \text{for all } m,\  n \leq N\delta, 
  \quad \text{for all } N \geq N_0 \,.
\end{equation}
Consider the martingale
\begin{equation*}
  M_n \deq e^{\bar{\theta} Y^N_n} 
  \left(\prod_{k = 0}^{n-1} \Phi^N( Y^N_k,k ) \right)^{-1} \,,
\end{equation*}
and set $R^N \deq  \inf\{n \geq 0 \,:\, Y^N_n \geq N\delta \}$. It
follows from \eqref{eq:phi_bound} that, on
the event $\{ T^N \leq R^N \wedge N\delta\}$, 
\[
M_{T^N} \geq 1\,, \quad \text{for all } N \geq N_0 \,.
\] 
Hence for $N \geq N_0$,
\[
e^{\theta} > \E M_0 = e^{\bar{\theta}} =
   \E M_{T^N \wedge R^N \wedge N\delta}
   \geq \Pb{T^N \leq R^N \wedge N\delta} \,.
\] 
However \eqref{eq:old_fluid_limit_result}  implies that, for 
$\delta < g/2$, $\Pb{R^N < T^N \leq N\delta} \rightarrow 0$, 
and \eqref{eq:unlikely_to_die_before_delta}  follows.
\end{proof}

\section{Identifiability In Patch-Free Processes}
\label{sec:patch-free_process}
We now focus on the case of patch-free hypergraph processes,
proving in this section Theorem \ref{thm:intro_no_patch_process}.

\subsection{A Coupled Family of Random Walks}
Let $\{P_t(n)\}_{t \geq  0}$, $n \in \N$, be a family of independent Poisson
processes, all of rate $2 \rho_2 > 0$, and consider the coupled family
of random walks $\{\xi_t(n)\}_{n \geq  0}$, for 
$t \in \R_+$, where $\xi_t(0) = 1$ for all $n$, and  
\begin{gather}
  \xi_t(n + 1) = \xi_t(n) + (P_t(n + 1) - 1) \one_{\{n < M_t\}} \,; \\ 
  M_t \deq \inf\{n \geq  0 \,:\,  \xi_t(n) = 0\}  
  \in \N \cup \{ \infty \} \,. \label{eq:rw_extinction_time} 
\end{gather}
 
The marginal law of $M_t$ is given by \eqref{eq:bp_population}
with $\beta_2 \deq t\rho_2$. There is a relation between
$\{\xi_t(n)\}_{n \geq  0}$ and the multigraph structure function:
since $g_2(t)$ is the largest root in $[0,1]$ of $2 t \rho_2 x +
\log(1 - x) = 0$, we have as a special case of
\eqref{eq:bp_population}: 

%\addtocounter{subsection}{1}
\begin{lem}
  The first time $t$ at which $\{\xi_t(n)\}_{n \geq  0}$ escapes to
  infinity is related to the multigraph lower envelope \eqref{eq:g2s}
  as follows:  
  \begin{equation}
    \Pb{M_t = \infty} = g_2(t) \,.
  \end{equation}
  Moreover $t \mapsto  M_t$ is an increasing process by the coupling,
  so $\chi \deq \inf \{t \geq  0 \,:\, M_t = \infty\}$ is a continuous
  random variable with distribution function $g_2(t)$.
\end{lem}

\subsection{Notation}
\label{sec:no_patch}
We finally turn to the case of a Poisson($\rho$) hypergraph process
$\{\Lambda_t\}_{t \geq 0}$ without patches, i.e.\ such that 
\[
\rho_0 \deq \rho_1 \deq 0 < \rho_2\,, 
\quad  \rho(x) \deq  \sum_{k \geq  2} \rho_k x^k \,, 
  x \in [0, 1] \,.
\] 
Write $T^N_t$ for the number of vertices in the domain of $v_0$ in
$\Lambda_t$, and write $Z^N_t$ for the number of hyperedges
identifiable from $v_0$ in $\Lambda_t$. Set  
$\bar{T}^N_t \deq N^{-1}T^N_t$ and 
$\bar{Z}^N_t \deq N^{-1}Z^N_t$. Using 
\eqref{eq:rw_extinction_time}, we define what will turn out to be the
macroscopic limits for Theorem \ref{thm:intro_no_patch_process}.
\begin{align*}
  \bar{T}_t & \deq g(t) \one_{\{M_t = \infty\}} \,; \\ 
  \bar{Z}_t & \deq \left\{ t\rho(g(t)) - [1 - g(t)]\log(1 - g(t))
  \right\} \one_{\{M_t = \infty\}} \,.
\end{align*}

\begin{proof}[Proof of Theorem \ref{thm:intro_no_patch_process}]
\hspace{0.1in}

\emph{Step I}.  Extending the notation of Theorem 
\ref{thm:identif_sans_patches} %\ref{thm:no_patch_process}
let $\Lambda_t(n)$ denote the hypergraph that
results from applying $n$ steps of randomized collapse to
$\Lambda_t + \one_{\{v_0\}}$; $Y^N_t(n)$ and $Z^N_t(n)$ count the number
of patches, and the amount of debris, respectively in
$\Lambda_t(n)$, and $n$ is assumed to satisfy: 
\[
n \leq  T^N_t \deq 
    \inf\{n \geq  0 \,:\,  Y^N_t(n) = 0 \} \,.
\] 
Consider a finite set of time points 
$0 < t_1 < \ldots < t_r$. 
The hypergraph collapses of  
$\Lambda_{t_1} + \one_{\{v_0\}}, \ldots,   
\Lambda_{t_r} + \one_{v_0}$ are coupled together as follows: perform
the $(n+1)$st step of randomized collapse by choosing a patch
uniformly at random from the smallest unstable hypergraph. Poisson
symmetries imply that this amounts to randomized collapse for each
of the unstable hypergraphs. Condition on the event: 
\begin{equation} \label{eq:condition_on_m_and_k}
  \bigcap_{i = 1}^{r} \{ Y^N_{t_i}(n) = m_i\,, 
  \ Z^N_{t_i}(n) = k_i \} \,.
\end{equation}
For $i$ such that $m_i = 0$, evidently $Y^N_{t_i}(n + 1) = 0$
and $Z^N_{t_i}(n + 1) = k_i$. 
For those $i$ such that $m_i \geq  1$, we may write:
\begin{align*}
  Y^N_{t_i}(n + 1) & = m_i - 1 - W^N_{t_i}(n + 1) + U^N_{t_i}(n + 1)\,; \\ 
  Z^N_{t_i}(n + 1) & = k_i + 1 + W^N_{t_i}(n + 1)\,,
\end{align*} 
where the random increments are distributed as follows.
Take $q$ to be the least $i \in \{1, 2, \ldots, r\}$ for which 
$m_i \geq 1$, and take $W^N_{t_q}(n + 1)$ and $U^N_{t_q}(n + 1)$
 independent such that 
\begin{equation} \begin{split}
    W^N_{t_q}(n + 1) & \sim {\rm Binomial}\left( m_q - 1,
    \frac{1}{N-n}\right) \,; \\
    U^N_{t_q}(n + 1) & \sim {\rm Poisson} \left( 
    (N - n - 1) t_q \lambda_2(N-n) \right) \,,
\end{split} \end{equation}
where $\lambda_2(N, n)$ is as in 
\eqref{eq:lambda2}. Because of the coupling, we may take subsequent
increments (for $i = q, \ldots, r - 1$) to be independent and of the form: 
\begin{align*}
  W^N_{t_{i+1}}(n + 1) - W^N_{t_i}(n + 1) & \sim
  {\rm Binomial} \left(m_{i+1} - m_i, \frac{1}{N-n} \right) \,; \\
  U^N_{t_{i+1}}(n + 1) - U^N_{t_i}(n + 1) & \sim
  {\rm Poisson} \left( (N - n - 1)(t_{i+1} - t_i) \lambda_2(N,n)
  \right) \,.
\end{align*}

\emph{Step II}. Observe that the behavior of $\lambda_2(N, n)$ depends
on whether $n \deq O(1)$, or $n \deq O(N)$. It follows from
\eqref{eq:lambda2_limit} and the calculations in Step I that, 
conditional on \eqref{eq:condition_on_m_and_k}, the joint law of 
\[
\left( (Y^N_{t_1}(n + 1), Z^N_{t_1}(n + 1)), \ldots, 
       (Y^N_{t_r}(n + 1), Z^N_{t_r}(n + 1)) \right)
\]
converges as $N \rightarrow \infty$ to the conditional law of
\[
\left( (\xi_{t_1}(n + 1), k_1 + 1), \ldots, 
    (\xi_{t_r}(n + 1), k_r + 1) \right)
\]
given that 
$\xi_{t_1}(n) = m_1, \ldots, \xi_{t_r}(n) = m_r$. 
Evidently $Z_{t_i}^N(0) = 0$ for all $i$.
Since $n$ was arbitrary, and since for each $t$ both 
$\{\xi_t(n)\}_{n \geq  0}$ and 
$\{ (Y^N_t(n), Z^N_t(n)) \}_{n \geq  0}$
 are Markov, we have now proved convergence in distribution as
$N \rightarrow \infty$: 
\begin{multline*}
  \{ (Y^N_{t_1}(n), Z^N_{t_1}(n)), \ldots, 
  (Y^N_{t_r}(n), Z^N_{t_r}(n)) \}_{n \geq  0} \\
  \weakconv
  \{ (\xi_{t_1}(n), n \wedge M_{t_1}), \ldots, 
  (\xi_{t_r}(n), n \wedge M_{t_r})\}_{n \geq  0} \,.
\end{multline*}
In particular, in the notation of 
\eqref{eq:rw_extinction_time} and Section
\ref{sec:no_patch},
\begin{equation} \label{eq:fd_distribs}
  \left( (T^N_{t_1}, Z^N_{t_1}), \ldots, 
  (T^N_{t_r}, Z^N_{t_r}) \right) 
  \weakconv
  \left( (M_{t_1}, M_{t_1}), \ldots, 
  (M_{t_r}, M_{t_r}) \right) \,.
\end{equation}

\emph{Step III}. To prove 
\eqref{eq:micro_no_patch_fluid_limit}
it suffices, in the light of 
\eqref{eq:fd_distribs},
to prove tightness of 
$\{(T^N_t, Z^N_t)\}_{t \geq 0}$ with
respect to the Skorohod topology of  
$D([0, \infty), (\N \cup \{\infty\})^2)$. 
On $(\N \cup \{\infty\})^2$, we shall use the metric
\[
d\left( (m,n),(p,q) \right) \deq
   \max\left\{\left|\frac{1}{p} - \frac{1}{m} \right|, 
   \left|\frac{1}{q} - \frac{1}{n} \right|\right\} \,
\]
understanding that $1/\infty = 0$. We shall 
verify the condition of Aldous for tightness of 
$\{(T^N_t, Z^N_t)\}_{t \geq  0}$, as stated 
in \ocite{B:CPM2}, p.\ 176, or \ocite{K:FMP}, 
p.\ 314, with respect to this metric. 
Since $s \mapsto T^N_s$ and $s \mapsto Z^N_s$ are
non-decreasing processes, the condition takes a slightly simpler form
than usual: it suffices to show that, for each $\epsilon > 0$ and
$\eta > 0$, there exist $h$ and $N_0$ such that for every bounded sequence
of optional times $\sigma^N$ with respect to 
$\{(T^N_t, Z^N_t)\}_{t \geq 0}$,
and for every $N \geq  N_0$, 
\begin{equation} \label{eq:aldous_crit_1}
  \Pb{\max\left\{\left|\frac{1}{T^N_{\sigma + h}}  - 
    \frac{1}{T^N_{\sigma}}\right|,
    \left|\frac{1}{Z^N_{\sigma + h}}  
    - \frac{1}{Z^N_{\sigma}}\right| \right\} \geq
    \ep } < \eta \,,
\end{equation}
where $\sigma$ is short for $\sigma^N$ in the subscripts.

Proposition \ref{prop:markov_property}
established
that $\{ (T^N_t, Z^N_t) \}_{t \geq  0}$ is a Markov process. By the strong
Markov property, the conditional law of 
$T^N_{\sigma + h} - m^N$, given
that $T^N_\sigma = m \deq m^N$, and $Z^N_\sigma = q^N$, is that
same as that of the number of identifiable vertices in a
Poisson($\hat{\beta}$) random hypergraph $\hat{\Lambda}^N$ 
on $\hat{N} \deq N - m$
vertices, where by the reasoning of Lemma 
\ref{lem:stopset} %\ref{sec:residual_hg},
and the fact that $\rho_1 = 0$,
\[
\hat{\beta}_1 \deq \frac{hN}{N - m}
    \sum_{k \geq 2} \rho_k \binom{m}{k - 1}
    \binom{N - m}{1} / \binom{N}{k} \,.
\]
Suppose $\epsilon > 0$ and $\eta > 0$ are given. In the case where 
$\min\{m^N, q^N\} > 1/\epsilon$, it follows that 
\begin{equation}
  \max\left\{ \left| \frac{1}{T^N_{\sigma + h}} 
  - \frac{1}{T^N_\sigma} \right|, \ 
  \left| \frac{1}{Z^N_{\sigma + h}}  - 
  \frac{1}{Z^N_\sigma} \right| \right\} < \epsilon \,.
\end{equation}

On the other hand, if $m^N \leq  1/\epsilon$, then 
\[
\hat{\beta}_1 \hat{N} \leq
   h N^2 \sum_{k \geq 2} \rho_k \binom{2/ \ep}{k -  1}
      / \binom{N}{k} = \frac{2 h \rho_2}{\ep} +O(N^{-1}) \,.
\]
Choose $N_0$ so large that, for $N \geq  N_0$, the right side is not more than 
$3 h \rho_2/\epsilon$; now it is true that, for any 
\[
h \leq \frac{-\epsilon \log(1 - \eta)}{3\rho_2} \,,
\]
and for any $N \geq  N_0$, the probability that $\hat{\Lambda}^N$ has no
patches, and hence no identifiable vertices nor identifiable hyperedges,
is at least $1-\eta$; in that case, $T^N_{\sigma + h}
= T^N_\sigma$ and 
and $Z^N_{\sigma + h} = Z^N_{\sigma}$. In summary, for such $N$ and
$h$,
\eqref{eq:aldous_crit_1} holds. Hence 
$\{(T^N_t, Z^N_t)\}_{t \geq 0}$ is tight, 
and \eqref{eq:micro_no_patch_fluid_limit}
follows. 

\emph{Step IV}. As for \eqref{eq:macro_no_patch_fluid_limit}
%\eqref{eq:micro_no_patch_fluid_limit},
we need only check the convergence of finite-dimensional
distributions, i.e.\ that 
\begin{equation} \label{eq:macro_fd_distribs}
  \left( (\bar{T}^N_{t_1}, \bar{Z}^N_{t_1}), \ldots, 
  (\bar{T}^N_{t_r}, \bar{Z}^N_{t_r}) \right)
  \weakconv
  \left( (\bar{T}_{t_1}, \bar{Z}_{t_1}), \ldots, 
  (\bar{T}_{t_r}, \bar{Z}_{t_r}) \right) \,. 
\end{equation}
for every finite set of time points $0 < t_1 < \ldots < t_r$. For the
case $r = 1$, the validity of \eqref{eq:macro_fd_distribs}
follows from Theorem 
\ref{thm:identif_sans_patches}.
For the sake of brevity, restrict our discussion of the
case $r > 1$ to the $\bar{T}$ component; the argument for 
the $\bar{Z}$ component is similar. It suffices to show,
for all $q = 2,\ldots,r$, and all
$\epsilon > 0$, that 
\begin{equation} \label{eq:jump_to_infinity}
  \P\!\left(\bigcap_{\substack{ i,j \\ 1 \leq i < q \leq  j \leq  r}} 
  \{ \bar{T}^N_{t_i} < \ep \} \cap
  \{ | \bar{T}^N_{t_j} - g(t_j)| < \ep \} \right)
\rightarrow
\Pb{M_{t_{q-1}} < \infty = M_{t_q}} \,.
\end{equation}

By our knowledge of the finite dimensional distributions from Theorem
\ref{thm:identif_sans_patches}, 
the left side of
is well approximated by 
\[
1 - \Pb{\bar{T}^N_{t_{q - 1}} \geq \epsilon}
- \Pb{\bar{T}^N_{t_q} \leq g(t_q) - \epsilon} \,,
\]
and for $\epsilon$ sufficiently small, this converges to the right
 side of \eqref{eq:jump_to_infinity}.
\end{proof} 

\section{Future Directions} \label{sec:future}

We have not explained here the role of the upper envelope
\eqref{eq:g_star_definition}, even though it was
included in the classification of structure functions. It is related
to dual hypergraph collapse and the size of the core, as in 
\ocite{C:core}. 
We shall give the corresponding asymptotic results in a
future paper. 

\section*{Acknowledgments} We thank Peter Matthews and the
referees for suggesting various expository improvements.

\end{document}